\documentclass[12pt]{amsart}
\usepackage{amssymb,amsmath,amsthm, color}
\theoremstyle{plain}
\newtheorem{thm}{Theorem}[section]
\newtheorem{cor}[thm]{Corollary}
\newtheorem{prop}[thm]{Proposition}
\newtheorem{lem}[thm]{Lemma}

\theoremstyle{definition}
\newtheorem{rem}[thm]{Remark}
\newtheorem{defn}[thm]{Definition}

\newtheorem{eg}[thm]{Example}

\renewcommand{\phi}{\varphi}
\newcommand{\lip}{\langle}
\newcommand{\rip}{\rangle}
\newcommand{\ip}[1]{\lip #1 \rip}

\newcommand{\ol}{\overline}

\newcommand{\upchi}{{\raise.35ex\hbox{$\chi$}}}
\newcommand{\sot}{\textsc{sot}}
\newcommand{\wot}{\textsc{wot}}
\newcommand{\AND}{\text{ and }}

\newcommand{\qand}{\quad\text{and}\quad}
\newcommand{\qor}{\quad\text{or}\quad}

\newcommand{\qforal}{\quad\text{for all}\quad}

\newcommand{\Dim}{\operatorname{dim}}
\newcommand{\diag}{\operatorname{diag}}

\newcommand{\ran}{\operatorname{Ran}}
\newcommand{\rank}{\operatorname{rank}}
\newcommand{\spn}{\operatorname{span}}
\newcommand{\Tr}{\operatorname{Tr}}
\newcommand{\bC}{{\mathbb{C}}}
\newcommand{\bM}{{\mathbb{M}}}

\newcommand{\bT}{{\mathbb{T}}}

\newcommand{\bZ}{{\mathbb{Z}}}
  \newcommand{\B}{{\mathcal{B}}}
  \newcommand{\C}{{\mathcal{C}}}
  \newcommand{\D}{{\mathcal{D}}}
  \newcommand{\E}{{\mathcal{E}}}
  \newcommand{\F}{{\mathcal{F}}}
  \newcommand{\G}{{\mathcal{G}}}
\renewcommand{\H}{{\mathcal{H}}}

  \newcommand{\K}{{\mathcal{K}}}  
\renewcommand{\L}{{\mathcal{L}}}
  \newcommand{\M}{{\mathcal{M}}}
  \newcommand{\N}{{\mathcal{N}}}

  \newcommand{\R}{{\mathcal{R}}}

  \newcommand{\U}{{\mathcal{U}}}
  
  \newcommand{\W}{{\mathcal{W}}}
  \newcommand{\X}{{\mathcal{X}}}
  \newcommand{\Y}{{\mathcal{Y}}}
  
\newcommand{\fD}{{\mathfrak{D}}}

\newcommand{\fH}{{\mathfrak{H}}}

\newcommand{\fK}{{\mathfrak{K}}}

\newcommand{\fT}{{\mathfrak{T}}}

\newcommand{\ep}{\varepsilon}

\begin{document}

\title{Transitive Spaces of Operators}

\author[K.R.Davidson]{Kenneth R. Davidson}
\address{Pure Math.\ Dept.\\U. Waterloo\\Waterloo, ON\;
N2L--3G1\\CANADA}
\email{krdavids@uwaterloo.ca \vspace{-2ex}}

\author[L.W.Marcoux]{Laurent W. Marcoux}
\email{lwmarcoux@uwaterloo.ca \vspace{-2ex}}

\author[H.Radjavi]{Heydar~Radjavi}
\email{hradjavi@uwaterloo.ca}

\thanks{Authors partially supported by NSERC grants.}
\thanks{\today}

\begin{abstract}
We investigate algebraic and topological transitivity and, more generally, $k$-transitivity for linear spaces of operators.  In finite dimensions, we determine minimal dimensions of $k$-transitive spaces for every $k$, and find relations between the degree of transitivity of a product or tensor product on the one hand and those of the factors on the other.   We present counterexamples to some natural conjectures.  Some infinite dimensional analogues are discussed.  A simple proof is given of Arveson's result on the weak-operator density of transitive spaces that are masa bimodules.
\end{abstract}

\date{}
\maketitle

A collection  $\E$   of operators  from a vector space $\X$ to another vector space $\Y$ 
is said to be \emph{transitive}  if, given nonzero $x \in \X$ and $y \in \Y$, 
there is an element $A \in \E$ such that $A x = y$.   
If $\E$ has some additional structure imposed on it, for example, if it is a ring or 
an algebra with $\X = \Y$,  then a great deal is implied by the transitivity assumption.   
This paper represents an effort to extract as much information as possible from 
transitivity (and its strengthening given below)  from the linear structure on $\E$ alone.

Let $\H$ and $\K$ be Hilbert spaces.        
A subspace $\L$ of $\B(\H,\K)$ is \textit{$k$-transitive} provided that 
for every choice of $k$ linearly independent vectors $x_1,\dots,x_k$ in $\H$ 
and arbitrary vectors $y_1,\dots,y_k$ in $\K$, 
there is an $A \in \M$ so that $Ax_i = y_i$ for $1 \le i \le k$.   Thus transitivity coincides with $1$-transitivity.  
$\L$ is \textit{topologically $k$-transitive} if this can be done approximately;  that is,
for each $\varepsilon > 0$ there exists $A \in \L$ such that $\| A x_i - y_i \| < \varepsilon$, 
$1 \leq i \leq k$.
In the finite dimensional setting, the two notions coincide.    As with transitivity, if $k =1$, we shall abbreviate the notation and refer to a space as being \textit{topologically transitive}.

Starting  with  finite dimensions, we present a number of positive results 
and a lot of counterexamples to natural conjectures.   
Among other things, we determine the minimal dimension of k-transitive 
subspaces of $\B(\H,\K)$ in terms of the dimensions of $\H$ and $\K$.   
We consider spans of products as well as tensor products of spaces, 
and study the relations between their  degree of transitivity and 
those of their constituent spaces.   
We investigate the relations between the minimal and maximal ranks 
present in a transitive subspace of $\bM_n$.
In particular, we show that such a subspace  must contain invertible elements.
In the infinite dimensional setting, there are fewer positive results 
and more counterexamples.   
We are able to extend some of the finite dimensional results.   
We provide a simple proof of a result of Arveson~\cite{Arv} 
that a topologically transitive subspace of $\B(\H,\K)$ which is a 
masa bimodule is \wot-dense in $\B(\H,\K)$.   

Any transitive subalgebra of $\bM_n=\B(\H_n)$, the space of 
$n \times n$ complex matrices, is equal to all of $\bM_n$ by Burnside's Theorem.
The situation is completely different for subspaces.  For every $0 \le k < \min \{m, n \}$,  there are $k$-transitive subspaces of $\bM_{mn}$ which are 
not\linebreak $(k+1)$-transitive (see Example~\ref{k_not_kplus1_trans}).   
In infinite dimensions, a topologically transitive operator algebra has
no proper invariant subspaces.  Whether it is \wot-dense in $\B(\H)$ is
the famous Transitive Algebra Problem, a generalization of the Invariant Subspace Problem.   
Again, if we consider subspaces,
there are many proper $k$-transitive \wot-closed subspaces.   

It is a well-known result of Azoff \cite{Az} that a subspace $\L$ of 
$\bM_{mn}=\B(\H_n,\H_m)$ is $k$-transitive if and only if the pre-annihilator 
$\L_\perp$ contains  no non-zero elements of rank at most $k$.
Here $(\bM_{mn})_*$ is identified with $\bM_{nm}$ equipped with the trace norm
via the bilinear pairing $\ip{A,T} = \Tr(AT)$.  
 
Indeed suppose that $0\ne T = \sum_{i=1}^k x_iy_i^*$ belongs to $\L_\perp$,
where the vectors $x_i$ are linearly independent.
Then the $k$-tuple $(Ax_1,\dots,Ax_k)$ is orthogonal to $(y_1,\dots,y_k)$
for all $A\in\L$. Thus $\L$ is not $k$-transitive.  
Conversely, if $\L$ is not $k$-transitive, then for some linearly
independent set $x_1,\dots,x_n$, the $k$-tuples $(Ax_1,\dots,Ax_k)$
span a proper subspace of $\H_m^k$, and thus is orthogonal to a
non-zero vector $(y_1,\dots,y_k)$.
Reversing the argument shows that $T = \sum_{i=1}^k x_iy_i^*$ belongs to
$\L_\perp$.

In the infinite-dimensional setting, we identify $(\B(\H, \K))_*$ with the space $\C_1(\K, \H)$ of trace class operators from $\K$ to $\H$ via the same bilinear pairing.    The above Theorem of Azoff then applies to topological $k$-transitivity, and shows that a subspace $\L$ of $\B(\H, \K)$ is topologically $k$-transitive if and only if the pre-annihilator $\L_\perp$ of $\L$ contains no non-zero elements of rank less than or equal to $k$.

We end this introduction with a simple but handy observation which will be used implicitly throughout this paper, namely:  if $\L \subseteq \B(\H, \K)$ is topologically $k$-transitive, and if $P \in \B(\H)$ and $Q \in \B(\K)$ are projections with $\min (\rank P, \rank Q) \geq k$, then $Q \L P \subseteq \B(P \H, Q \K)$ is also topologically $k$-transitive.   (This follows immediately from Azoff's Theorem.)   If, furthermore, $P$ and $Q$ are finite rank, then $Q \L P$ is in fact $k$-transitive.

\section{Dimension}

In this section, we find the minimal dimension of a $k$-transitive subspace of $\bM_{mn}$.  
By Azoff's Theorem as outlined above, if $k \geq \min (m, n)$, 
then $\L_\perp = 0$ and hence $\L = \bM_{mn}$.   
As such, we always assume that $k < \min (m,n)$.

\begin{lem}\label{diag}
For every $p\ge 1$ and every $0 \le k \le p$, the diagonal subalgebra $\fD_p$ of $\bM_p$ contains a subspace of dimension $p-k$
containing no non-zero element of rank at most $k$.
\end{lem}

\begin{proof} It suffices to choose $p-k$ diagonal operators with the property
that when restricted to any $p-k$ diagonal entries, they are linearly independent.
For then a linear combination which has $p-k$ zeros must be zero.
An example of such a sequence is $D_j = \mathrm{diag}(1^j,2^j,\dots,p^j)$ for $0 \le j < p-k$.
\end{proof}

The special case of the following result for $k=1$ was established by Azoff \cite{Az}.

\begin{thm}\label{transitive}
The minimal dimension of a $k$-transitive subspace of $\bM_{mn}$
is $k(m+n-k)$.
\end{thm}

\begin{proof}
Using Azoff's Theorem, we search instead for
the maximal dimension of the pre-annihilator of a $k$-transitive subspace $\L$.
This subspace $\L_\perp$ cannot intersect the closed variety $\R_k$ of
matrices of rank at most $k$ except in $\{0\}$.
This variety has dimension $k(m+n-k)$ (see \cite[Prop.~12.2]{JH}).
It follows that $\dim \L_\perp + \dim \R_k \le mn$.
Hence $\dim \L \ge \dim \R_k = k(m+n-k)$.

On the other hand, consider the subspace $\N$ of $\bM_{n m}$ which has zeros on diagonals
of length $p \le k$ and has dimension $p-k$ on the diagonals with 
$p > k$ entries such that the rank of any non-zero element on one such
diagonal is always at least $k+1$.  
This is possible by Lemma~\ref{transitive}.
There are $m+n-1$ diagonals, and $2k$ of them have length at most $k$.
So the dimension of $\N$ is
\begin{align*}
  mn - (m \!+\! n \!-\! 1 \!-\! 2k)k - 2 \sum_{p=1}^{k} p  &=
 mn - (m \!+\! n \!-\! 1 \!-\! 2k)k - (k \!+\! 1)k \\ &= mn - k(m + n - k)
\end{align*}
Thus $\L = \N^\perp$, the annihilator of $\N$,  has dimension $k(m+n-k)$.

Consider a non-zero element $N$ of $\N$.
It must be non-zero on some diagonal.  
Let $p_0$ be the shortest non-zero diagonal.
Consider the square submatrix containing the $p_0$th diagonal
as its main diagonal.  
This submatrix is triangular, and hence its rank is at least
as great as the rank of the diagonal, which is at least $k+1$.
Hence $\N$ contains no non-zero elements of rank at most $k$.
Therefore $\L$ is $k$-transitive.
\end{proof}

\begin{eg} \label{trace_zero_matrices}
By the above theorem, any $(n-1)$-transitive subspace $\L$ of $\bM_n$ 
necessarily has dimension $n^2 - 1$.   
It is then not hard to see that we can find invertible elements 
$S, T \in \bM_n$ so that $S \L T = \mathfrak{sl}_n$, 
the space of trace zero matrices in $\bM_n$.

Since multiplying a $k$-transitive subspace $\M$ of $\B(\H, \K)$ 
by invertibles $S \in \B(\K, \K^\prime)$ and $T \in \B(\H^\prime, \H)$ 
yields a $k$-transitive subspace of $\B(\H^\prime, \K^\prime)$, 
we shall think of $\M$ and $S \M T$ as being \emph{equivalent} 
insofar as transitivity is concerned.  

A similar statement holds for topologically $k$-transitive spaces.
\end{eg}

\begin{eg} \label{k_not_kplus1_trans}
As mentioned in the introduction, there are $k$-transitive subspaces of 
$\B(\H, \K)$ which are not $(k+1)$-transitive provided that 
$k < \min \{ \Dim \H,\Dim \K \}$.    
For example, if $\H = \H_n$ and $\K = \H_m$, then it suffices to 
consider $R \in \bM_{nm}$ of rank $k+1$.   
With $\M = \spn \{ R \}$, $\L = \M^\perp \subseteq \bM_{mn}$ is such an example.  

This example can easily be adapted to the infinite-dimensional setting to produce subspaces $\L \subseteq \B(\H, \K)$ which are topologically $k$-transitive but not topologically $(k+1)$-transitive.
\end{eg}

\begin{eg} \label{toeplitz_are_trans}
The space $\mathfrak{T}_n$ of all Toeplitz matrices $T=\big[t_{i-j}\big]$ in $\bM_n$ 
is a transitive subspace of dimension $2n-1$.   It is routine to verify that the pre-annihilator of $\mathfrak{T}_n$ consists of those matrices in $\bM_n$ whose entries along any diagonal sum to zero.   The rank of such a matrix is at least as big as the length of the smallest non-zero diagonal, which is at least two if the matrix is non-zero.   

By Theorem~\ref{transitive}, this is the minimal possible dimension of a transitive subspace of $\bM_n$.
\end{eg}

\section{Dually Transitive}

In finite dimensions, we can consider both $\L$ and $\L_\perp$ as spaces of matrices.
So we can ask whether both can be $k$-transitive.  
Dimension arguments show that this requires the space to be sufficiently large.

\begin{prop}\label{dual_trans}
There is a subspace $\L$ of $\bM_4$ such that both $\L$ and $\L_\perp$
are transitive.
\end{prop}

\begin{proof} 
It suffices to exhibit a transitive subspace $\L$ so that neither $\L$ nor
$\L_\perp$ contains a rank one element. 
Let $\Phi$ be a linear bijection of $\bM_2$ onto itself with four
distinct eigenvalues whose corresponding eigenvectors all have rank 2.
Define
\[ \L = \Big\{ \begin{bmatrix}A&\Phi(B)\\ B&A\end{bmatrix} : A,B \in \bM_2 \Big\} .\]

Suppose that $\L$ contains a rank 1 element 
$\begin{bmatrix}xu^*&xv^*\\ yu^*&yv^*\end{bmatrix}$,
where $x$, $y$, $u$, $v$ are vectors in $\bC^2$.
Then comparing the diagonal entries shows that $xu^*=yv^*$,
so that $xv^*$ is a multiple of $yu^*$, say $xv^* = \lambda yu^*$.
So now a comparison of the off-diagonal entries shows that
$\Phi(yu^*) = \lambda yu^*$, contrary to fact.
Therefore $\L$ contains no rank 1 elements.

The analysis of $\L_\perp$ is similar since $\Phi^t$ has the same
eigenvector property.
\end{proof}

\begin{eg}\label{Ex:dual_trans}
For example, let $\Phi\Big(\begin{bmatrix}a&b\\ c&d\end{bmatrix}\Big) = 
\begin{bmatrix}d&2c\\ b&a\end{bmatrix}$.
Then 
\[ \L =  \begin{bmatrix}a&b&h&2g\\c&d&f&e\\e&f&a&b\\g&h&c&d\end{bmatrix}
 \qand
 \L_\perp =
 \begin{bmatrix}a&b&-h&-g\\c&d&-f&-e\\e&f&-a&-b\\ g/2 &h&-c&-d\end{bmatrix}
\]
\end{eg}

\smallskip

Here is a general technique modelled on the previous example.

\begin{thm}
If  $1 \le k < n/(2+ \sqrt2)$, then $\bM_{2n}$ contains a subspace $\L$
such that both $\L$ and $\L_\perp$ are $k$-transitive.
\end{thm}

\begin{proof}
The proof of Theorem~\ref{transitive} shows that there is a subspace $\N$ of
$\bM_n$ with dimension $(n-k)^2$ which contains no non-zero elements of rank $k$ or less. 
Since $n >(2+ \sqrt2)k$, it is easy to check that $(n-k)^2 \ge n^2/2$.
Thus there is an injective linear map $T: \bM_n/\N \to \N$.
Let $\Phi = JTQ$, where $Q$ is the quotient map of $\bM_n$ onto $\bM_n/\N$
and $J$ is the injection of $\N$ into $\bM_n$.

Observe that for $A \in \bM_n$,
\[ \rank A \ge k+1 \qor \rank \Phi(A) \ge k+1 .\]
Indeed, if $\Phi(A) \ne 0$, then it has rank at least $k+1$.
If $\Phi(A) = 0$, then $A \in \N$; so it has rank at least $k+1$.

Define a subspace $\L$ of $\bM_{2n}$ consisting of all elements of the form
$\begin{bmatrix}A&B\\\Phi(B)&\Phi(A)\end{bmatrix}$ 
where $A,B \in \bM_n$ are arbitrary.
A non-zero element of $\L$ has either $A$ or $B$ non-zero.
So at least one of the four matrix entries has rank at least $k+1$.
In particular, $\L_\perp$ is $k$ transitive.

Note that $\L_\perp$ consists of all matrices of the form
$\begin{bmatrix}\Phi^t(X)&-Y \\ \Phi^t(Y)&-X\end{bmatrix}$ 
where $X,Y \in \bM_n$ are arbitrary.
So the same argument shows that $\L_\perp$ contains no
non-zero elements of rank at most $k$.
Therefore $\L$ is $k$-transitive.
\end{proof}

\section{Tensor Products}

In this section, we consider tensor products of $k$-transitive subspaces.
The first lemma is well known, but is included for the convenience of the reader.

\begin{lem}\label{L:dualtensor}
Let $\L$ and $\M$ be subspaces of $\bM_{lp}$ and $\bM_{mn}$ respectively.   
Then $(\L \otimes \M)_\perp = \L_\perp \otimes \bM_{nm} + \bM_{pl} \otimes \M_\perp$.
\end{lem}

\begin{proof}
Clearly  $ \L_\perp \otimes \bM_{nm} + \bM_{pl} \otimes \M_\perp$ 
is contained in $(\L \otimes \M)_\perp$.   
The other containment follows from a simple dimension argument.
Indeed, let $\Dim \L = d_1$ and $\Dim \M = d_2$.   
Then $\Dim \L \otimes \M = d_1 d_2$, from which we deduce that 
\[ \Dim (\L \otimes \M)_\perp = lpmn - d_1 d_2 .\]
Now 
$
 \Dim (\L_\perp \otimes \bM_{nm}) = (lp - d_1) mn$ and $
 \Dim (\bM_{pl} \otimes \M_\perp) = lp (mn - d_2) .  
$
Since 
\[
 \L_\perp \otimes \bM_{nm} \cap \bM_{pl} \otimes \M_\perp 
 = \L_\perp \otimes \M_\perp
\]
has dimension $(lp - d_1)(mn - d_2)$, it follows that 
\begin{align*}
\Dim (\L_\perp \otimes \bM_{nm} + \bM_{pl} \otimes \M_\perp) 
	&= (lp - d_1) mn + lp (mn - d_2)\\
	& \ \ \ \  - (lp - d_1)(mn - d_2) \\
	&= lpmn - d_1 d_2. 
\end{align*}
 From this the result follows easily.	
\end{proof}

The main theorem of this section shows that under additional hypotheses, 
tensoring preserves some level of transitivity.

\begin{thm}\label{T:tensor}
Suppose that a subspace $\L \subset \bM_{lp}$ is $k$-transitive 
and is spanned by its rank $r$ elements.
If a subspace $\M \subset \bM_{mn}$ is $rk$-transitive, 
then $\L \otimes \M$ is $k$-transitive.
\end{thm}

\begin{proof}
By Lemma~\ref{L:dualtensor}, 
$(\L \otimes \M)_\perp = \L_\perp \otimes \bM_{nm} + \bM_{pl} \otimes \M_\perp$.
Consider $T \in (\L \otimes \M)_\perp$ as a $p \times l$ matrix $T = \big[ T_{ji} \big]$
with coefficients in $\bM_{nm}$.
Thus we may decompose $T_{ji} = R_{ji}+S_{ji}$ 
where $R = \big[ R_{ji} \big] \in \L_\perp \otimes \bM_{nm}$
and $S = \big[ S_{ji} \big] \in \bM_{pl} \otimes \M_\perp$.
The latter condition just says that each $S_{ji}$ is in $\M_\perp$. 
One detects that $R$ is in $\L_\perp \otimes \bM_{nm}$ by the fact that
it satisfies
\[  \sum_{i = 1}^p \sum_{j=1}^l l_{ij}  R_{ji} = 0 \qforal L = \big[ l_{ij} \big] \in \L .\]

For $L = \big[ l_{i j} \big] \in \L$, consider 
\begin{align}\label{E1}
 \sum_{i = 1}^p \sum_{j=1}^l l_{ij}  T_{ji} = 
 \sum_{i = 1}^p \sum_{j=1}^l l_{ij}  (R_{ji} + S_{ji}) = 
 \sum_{i = 1}^p \sum_{j=1}^l l_{ij}  S_{ ji} \in \M_\perp .
\end{align}
Suppose that $\rank T= 1$.
Then there are vectors $u_j$ and $v_i$ for $1 \le j \le p$ and $1 \le i \le l$
so that $T_{ji} = u_jv_i^*$.
Now if $L$ is rank one, then there are scalars $x_i$ and $y_j$ so that
$l_{ij} = x_iy_j$.  In this case,
\[
 \sum_{i = 1}^p \sum_{j=1}^l l_{ij}  T_{ji} = 
 \sum_{i = 1}^p \sum_{j=1}^l x_iy_j u_jv_i^* =
 \big(\sum_{j=1}^l y_j u_j\big) \big( \sum_{i = 1}^p \overline{x_i} v_i \big)^*
\]
is a rank one matrix.  It is then easy to see that if $\rank T \le k$
and $\rank L \le r$, then this sum has rank at most $rk$.

Clearly it suffices to satisfy equation (\ref{E1}) for a spanning subset of $\L$.
So we may suppose that each $L$ has rank at most $r$.
Consider an element $T \in (\L\otimes\M)_\perp$ with rank at most $k$.
The analysis of the previous paragraph yields an element of $\M_\perp$
with rank at most $rk$.
As $\M$ is $rk$-transitive, this means that these sums are all zero:
\[ \sum_{i = 1}^p \sum_{j=1}^l l_{ij}  S_{ ji} = 0 \qforal L = \big[ l_{i j} \big] \in \L .\]
Consequently, $S \in \L_\perp \otimes \bM_{nm}$.   Since $R$ is also in this set,
we conclude that  $T\in  \L_\perp \otimes \bM_{nm}$.

But $\L_\perp$ admits no non-zero operators of rank at most $k$, 
and thus neither does $\L_\perp \otimes \bM_{nm}$.  
To see this, think of such matrices as $n \times m$ 
matrices with coefficients in $\L_\perp$.
Any non-zero coefficient results in rank at least $k+1$.
Hence $(\L \otimes \M)_\perp$ contains no non-zero elements 
of rank at most $k$, and therefore $\L \otimes \M$ is $k$-transitive.
\end{proof}

\begin{defn}
A subspace $\L \subset \bM_{mn}$ is \textit{fully $k$-transitive} if
$\L \otimes \M$ is $k$-transitive whenever $\M$ is $k$-transitive.
\end{defn}

The following corollary yields a large class of fully $k$-transitive subspaces.
This is an immediate application of Theorem~\ref{T:tensor} taking $r=1$.

\begin{cor}\label{C:fullyktrans}
If $\L\subset \bM_{lp}$ is a $k$-transitive subspace which is spanned by
its rank one elements, then $\L$ is fully $k$-transitive.
\end{cor}
\begin{cor}
If $\L$ and $\M$ are both $k$-transitive spaces spanned by their
rank one elements, then $\L\otimes\M$ is fully $k$-transitive.
\end{cor}

\begin{proof}
$\L\otimes\M$ is $k$-transitive by Corollary~\ref{C:fullyktrans}.
The tensor product of rank one elements is rank one, and so
$\L\otimes\M$ is spanned by rank ones.
Hence it is fully $k$-transitive by the same corollary.
\end{proof}

Another easy consequence uses the fact that $\L \subset \bM_{lp}$
is always spanned by elements of rank at most $\min\{l,p\}$.

\begin{cor}
If $\L\subset \bM_{lp}$ is a $k$-transitive subspace 
and $\M \subset \bM_{mn}$ is $\min\{kl,kp\}$-transitive,
then $\L \otimes \M$ is $k$-transitive.
\end{cor}

Fully $k$-transitive spaces have a certain permanence.

\begin{prop}
If  $\L \subset \bM_{lp}$ is fully $k$-transitive, 
and $P,Q$ are idempotents in $\bM_l$ and $\bM_p$ respectively,
then $P\L Q \subset \B(Q\H_p,P\H_l)$ is fully $k$-transitive.
\end{prop}

\begin{proof} Let $\M \subset \bM_{mn}$ be a $k$-transitive space.
If $A \in (P\L Q)_\perp$, then $0 = \Tr(PLQA) = \Tr(LQAP)$ for all $L \in \L$.
So $QAP \in \L_\perp$, where $QAP$ is just $A$ considered as an element
of $\bM_{lp}$ rather than $\B(P\H_l,Q\H_p)$.  Therefore
\begin{align*}
 (P\L Q \otimes \M)_\perp &= 
 (P\L Q)_\perp \otimes \bM_{nm} + Q\bM_{pl}P \otimes \M_\perp \\
 &\subset \L_\perp \otimes \bM_{nm} + \bM_{pl} \otimes \M_\perp 
 = (\L \otimes \M)_\perp 
\end{align*}
The right hand side contains no rank $k$ matrices, and so neither does the left side.
Therefore $P\L Q$ is fully $k$-transitive.
\end{proof}

\begin{eg}
The space of Toeplitz matrices $\mathfrak{T}_m$ is fully transitive because
the rank one matrices $\big[ a^{i-j} \big]$ for $a \in \bC$ span $\mathfrak{T}_m$.
To see this, just observe that the entries on the first row and column
determine $T$, and that this may be any vector in $\bC^{2m-1}$.
The rank one matrices mentioned above correspond to the vectors 
$(a^k)_{ |k|<m}$.  Any choice of $2m-1$ distinct non-zero values of $a$
yields a basis.
\end{eg}

\begin{eg}
The space $\mathfrak{sl}_d \subset \bM_d$ from Example~\ref{trace_zero_matrices} is fully $(d-1)$-transitive since it is spanned by the rank one matrices 
$\{ E_{ij} : i \ne j \}$ and $\{E_{11}+E_{1j}-E_{j1}-E_{jj}: 2 \le j \le m\}$.
\end{eg}

\begin{eg}
Consider $\L = \mathfrak{sl}_d \otimes \bM_m \subset \bM_{dm}$.
We may think of this space as those $d\times d$ matrices $A=\big[ A_{ij} \big]$
with coefficients in $\bM_m$ such that $\sum_{i=1}^d A_{ii} = 0$.
Then $\L_\perp = (\mathfrak{sl}_d)_\perp \otimes \bM_m = \bC I_d \otimes \bM_m$.
Thus the minimum rank of a non-zero element of $\L_\perp$ is $d$.
Therefore $\L$ is $(d-1)$-transitive.
It is spanned by its rank one elements since both $\mathfrak{sl}_d$ and $\bM_m$ are.
Hence $\L$ is fully $(d-1)$-transitive.
\end{eg}

\begin{eg} There are transitive spaces which are not fully transitive.
Consider the space $\L$ from Example~\ref{Ex:dual_trans}. 
Evidently, the smallest rank of a non-zero element of $\L$ is 2.
We will show that $(\L_\perp \otimes \L_\perp)_\perp$ contains a rank one.
By symmetry, it follows that neither $\L \otimes \L$ 
nor $\L_\perp \otimes \L_\perp$ is transitive.
By Lemma~\ref{L:dualtensor}, 
\[ (\L_\perp \otimes \L_\perp)_\perp = \L \otimes \bM_4 + \bM_4 \otimes \L .\]

To find a rank 1 in this space, we look for matrices
$A,\dots,H$ in $\bM_4$, $X_{ij} \in \L$ and 
vectors $u_i,v_j \in \bC^4$  for $1 \le i,j \le 4$ so that
\[
 \begin{bmatrix}A&B&H&2G\\C&D&F&E\\E&F&A&B\\G&H&C&D\end{bmatrix} =
 \begin{bmatrix} & & &  \\ &X_{ij} + u_iv_j^*& & \\ & & & \\ & & & \end{bmatrix}
\]
This is equivalent to solving the system
\begin{alignat*}{2}
 u_1v_1^* - u_3v_3^* &= X_{33}-X_{11} \in \L 
    &\qquad u_2v_4^* - \phantom{2}u_3v_1^* &= \phantom{2}X_{31}-X_{24} \in \L\\
 u_1v_2^* - u_3v_4^* &= X_{34}-X_{12} \in \L 
    &\qquad u_2v_3^* - \phantom{2}u_3v_2^* &= \phantom{2}X_{32}-X_{23} \in \L\\
 u_2v_1^* - u_4v_3^* &= X_{43}-X_{21} \in \L 
    &\qquad u_1v_3^* - \phantom{2}u_4v_2^* &= \phantom{2}X_{42}-X_{13} \in \L\\
 u_2v_2^* - u_4v_4^* &= X_{44}-X_{22} \in \L 
    &\qquad u_1v_4^* - 2u_4v_1^* &= 2X_{41}-X_{14} \in \L
\end{alignat*}

Let $e_1,\dots,e_4$ be the standard basis for $\bC^4$.
One can check that 
\[
 \begin{bmatrix}u_1\\u_2\\u_3\\u_4\end{bmatrix} =
 \begin{bmatrix}e_4\\e_3\\e_2\\e_1\end{bmatrix} \qand
 \begin{bmatrix}v_1\\v_2\\v_3\\v_4\end{bmatrix} =
 \begin{bmatrix}\phantom{-}e_4\\\phantom{-}e_3\\-e_2\\-e_1\end{bmatrix}
\]
is a non-trivial solution.
\end{eg}

\section{Spans of products} \label{section4}

While transitive subspaces are not algebras, an algebra can be obtained by 
taking spans of products.
Thus in the matrix case, one eventually obtains $\bM_n$.
In many cases, this happens very quickly.
In particular, the order of transitivity increases quickly.

\begin{prop}\label{prod_span1}
If $\L \subset \bM_l$ is a transitive subspace which is spanned by its rank one elements,
then $\spn \L^2 = \bM_l$.
\end{prop}

\begin{proof}
$\L$ has no kernel; so we may select $l$ non-zero rank one 
elements $R_i = x_i y_i^*$ in $\L$ such that
$\{y_1,\dots,y_l\}$ forms a basis for $\H_l$.
Every matrix in $\bM_l$ may be written as $T=\sum_{i=1}^l u_i y_i^*$.
Choose $A_i \in \L$ so that $A_i x_i = u_i$.
Then $T = \sum_{i=1}^l A_i R_i$ belongs to $ \spn \L^2$.
\end{proof}

\begin{eg}\label{prodM}
Consider the space $\L$ from Example~\ref{Ex:dual_trans}. 
Let 
\[
 X = \begin{bmatrix}a_1&b_1&h_1&2g_1\\c_1&d_1&f_1&e_1\\
 e_1&f_1&a_1&b_1\\g_1&h_1&c_1&d_1\end{bmatrix}
 \qand 
 Y =   \begin{bmatrix}a_2&b_2&h_2&2g_2\\c_2&d_2&f_2&e_2\\
 e_2&f_2&a_2&b_2\\g_2&h_2&c_2&d_2\end{bmatrix}
\]
be arbitrary elements of $\L$.
Let $Z = X Y$, and consider the diagonal elements $z_{ii}$ of $Z$, $1 \leq i \leq 4$.
Then 
\begin{align*}
z_{11} &=  a_1 a_2 + b_1 c_2 + h_1 e_2 + 2 g_1 g_2 \\
z_{22} &=  c_1 b_2 + d_1 d_2 + f_1 f_2 + e_1 h_2 \\
z_{33} &=  e_1 h_2 + f_1 f_2 + a_1 a_2 + b_1 c_2 \\
z_{44} &=  2 g_1 g_2 + h_1 e_2 + c_1 b_2 + d_1 d_2.
\end{align*}

It follows that the expectation $E$ from the space $\spn \{\L, \L^2 \}$ onto 
the diagonal has range spanned by the three diagonal matrices 
\[
 D_1 =\diag(1, 0, 1, 0),\ 
 D_2 =\diag (1, 0, 0, 1) \AND
 D_3 =\diag(0, 1, 0, 1),
\] 
which is only three dimensional.    Indeed, with $X$ and $Z$ as above, $E(X)$ is spanned by $D_1$ and $D_3$, while 
\begin{align*}
E(Z) &= (a_1 a_2 + b_1 c_2) D_1 + (h_1 e_2 + 2 g_1 g_2) D_2 \\
	&\quad  + (c_1 b_2 + d_1 d_2) D_3 + (f_1 f_2 + e_1 h_2) (D_1 + D_3 - D_2). 
\end{align*}
In particular, $\spn \{ \L, \L^2 \} \not = \bM_4$.
Note that $\spn \{ \L, \L^2, \L^3\} = \bM_4$.
\end{eg}

The following concept is a substantial weakening of the notion of $k$-transitivity.

\begin{defn}
A subspace $\L \subset \B(\H, \K)$ is \textit{$k$-separating} if for every 
set $x_1,\dots,x_{k}$ of independent vectors in $\H$, there is an $L \in \L$
so that $Lx_i = 0$ for $1 \le i \le k-1$ and $Lx_{k} \ne 0$.
\end{defn}

The simple result  below shows why the property of being $k$-separating 
is nice to have when considering products of spaces.

\begin{lem} \label{easy_lem}
Suppose that $\L_1, \L_2$ are subspaces of $\bM_n$ with 
$\L_1$ transitive and $L_2$ $k$-separating.   
Then $\spn \L_1 \L_2 $ is $k$-transitive.
\end{lem}

\begin{proof}
Let $x_1,\dots,x_k$ be independent vectors, and let vectors $y_1,\dots,y_k$ be given.
Use the $k$-separating property to select elements 
$B_1,\dots,B_k$ in $\L_2$ such that $B_ix_j = \delta_{ij}z_i$,
where $z_i$ are non-zero vectors.
By the transitivity of $\L_1$, select elements $A_1,\dots,A_k$ in $\L_1$
so that $A_i z_i = y_i$.
Then $\sum_{i=1}^k A_iB_i$ belongs to $\spn \L_1 \L_2$ and takes
$x_i$ to $y_i$ for $1 \le i \le k$.
Therefore this space is $k$-transitive.
\end{proof}

Since $\bM_n$ is the unique $n$-transitive subspace of itself, we obtain:

\begin{cor}\label{easy_cor}
Suppose that $\L_1, \L_2$ are subspaces of $\bM_n$ with 
$\L_1$ transitive and $L_2$ $n$-separating.   
Then $\spn \L_1 \L_2 =\bM_n$.
\end{cor}

\begin{lem}\label{trans_implies_separating}
Let $1 \le k < \min\{m,n\}$.
Then a $k$-transitive subspace $\L$ of $\bM_{mn}$ is $(k+1)$-separating.
\end{lem}

\begin{proof}
Let $x_1,\dots,x_{k+1}$ be linearly independent in $\bC^n$; 
and set $\X = \spn\{x_1,\dots,x_{k+1}\}$.
Then the restriction $\M = \L|_\X \subset \bM_{m,k+1}$ is $k$-transitive.
By Theorem~\ref{transitive}, $\Dim\M \ge k(m+1)$.
The subspace of $\M$ which vanishes on $\spn\{x_1,\dots,x_k\}$
has $km$ linear conditions imposed, and hence it has dimension
at least $k$.  Thus there are elements of $\M$ which annihilate
$x_1,\dots,x_k$ and are non-zero on $x_{k+1}$.
\end{proof}


\begin{eg}\label{Toeplitz2sepnot3}
The set $\fT_n$ of $n \times n$ Toeplitz matrices ($n \geq 3$) is 
an example of a space which is $2$-separating but not $3$-separating.   
Indeed, the fact that $\fT_n$ is $1$-transitive implies that 
it is $2$-separating by the above result.  
On the other hand, if the first and last columns of a matrix 
$T$ in $\fT_n$ are both zero, then all entries of $T$ are zero, 
so that $\fT_n$ is not $3$-separating.
\end{eg}

\begin{eg}\label{nsepnottrans}
The converse of Lemma~\ref{trans_implies_separating} is false.

If $\M$ is a $k$-transitive but not $(k+1)$-transitive subspace
of $\bM_{mn}$, consider the subspace $\L$ of $\bM_{m+1,n}$ of the form
$L = \begin{bmatrix}x\\ M\end{bmatrix}$ where $x$ is an arbitrary row
vector in $\bC^n$ and $M \in \M$.
Then this is an $n$-separating space which is $k$-transitive but not
$(k+1)$-transitive.

In particular, if we take $\M=\{0\}$, 
then $\L$ is $n$-separating but is not even $1$-transitive.
\end{eg}


\begin{thm}\label{prod_increases_trans}
Suppose that subspaces $\K \subset \bM_{mn}$ and $\L \subset \bM_{np}$
are $k$-transitive and $l$-transitive respectively.
Then the product $\spn \K\L$ is $\min\{k+l,m,p\}$-transitive.
\end{thm}

\begin{proof}
We know that $\K_\perp$ contains no non-zero element of rank at most $k$;
and $\L_\perp$ contains no non-zero element of rank at most $l$.

Assume first that $l < \min\{n,p\}$.
Suppose that $A \in (\K\L)_\perp$ satisfies $1 \le \rank A \le k+l$.
Then $0 = \Tr(KLA)$ for all $K \in\K$ and $L\in\L$.
Hence $LA \in \K_\perp$ for all $L \in \L$.
As $\L$ is $(l+1)$-separating by Lemma~\ref{trans_implies_separating},
select $L \in \L$ which is non-zero on some vector in the range of $A$ and
annihilates $\min\{l,\rank A - 1 \}$ independent vectors in the range of $A$.
Then 
\[ 1 \le \rank LA \le (k+l) - l = k .\]
This contradicts the $k$-transitivity of $K$.
Thus $\spn\K\L$ is $(k+l)$-transitive if $k+l < \min\{m,p\}$.
But if $k+l \ge \min\{m,p\}$, then this shows that $(\K\L)_\perp = \{0\}$.
Hence $\spn\K\L = \bM_{mp}$ is $\min\{m,p\}$-transitive.

We obtain a similar conclusion if $k< \min\{m,n\}$.
If $k=\min\{m,n\}$ and  $l = \min\{n,p\}$, then
$\K=\bM_{mn}$ and $\L=\bM_{np}$.
Thus $\spn \K\L = \bM_{mp}$ is $\min\{m,p\}$ transitive.
\end{proof}

\smallskip
For $\M$ a subspace of $\bM_n$, let $\M^* = \{ M^*: M \in \M \}$, again considered as a subspace of  $\bM_n$.

\begin{prop} \label{Proposition4.10}
If $\L \subset \bM_n$ is transitive and spanned by its rank $r$ elements
and $\M^* \subset \bM_n$ is $r$-separating, then $\spn \L \M = \bM_n$.
\end{prop}

\begin{proof}
If $L\in\L$ has rank at most $r$, then for any vector $u \in L\H$ and $0 \ne x \in \H$,
we will show that $ux^* \in \L\M$.
Indeed, we may write $L = \sum_{i=1}^s u_i v_i^*$ where $u_1=u$,
$s \le r$ and $\{v_1,\dots,v_s\}$ are linearly independent.
Select $M\in\M$ so that $M^*v_1 = x$ and $M^*v_i = 0$ for $2 \le i \le s$.
Then 
\[ LM = \sum_{i=1}^s u_i v_i^* M = \sum_{i=1}^s u_i (M^*v_i)^* = ux^* .\]
As the ranges of elements (of rank at most $r$) of $\L$ span $\H$, the result follows.
\end{proof}


If $\M \subseteq \bM_n$ is $r$-transitive for some $r \geq 1$, then so is 
$\M ^*$.
This is easily seen by considering Azoff's characterisation of $r$-transitivity 
in terms of the preannihilator of $\M$.  
It then follows from Lemma~\ref{trans_implies_separating} that $\M^*$ 
is $(r+1)$-separating.
When $r=0$, no such statement holds.  So we obtain:

\begin{cor} 
If $\L \subset \bM_n$ is transitive and spanned by its rank $r$ elements
and $\M \subset \bM_n$ is $\max\{r-1,1\}$-transitive, then $\spn \L \M = \bM_n$.
\end{cor}

This allows an extension of Proposition~\ref{prod_span1} and Example~\ref{prodM}.

\begin{cor}
If $\L \subset \bM_n$ is transitive and is spanned by its rank $r$ elements,
then $\spn \L^{r+1} = \bM_n$.
\end{cor}

\begin{proof}
By Theorem~\ref{prod_increases_trans}, $\spn\L^r$ is $r$-transitive.
So by the previous lemma, $\spn \L^{r+1} = \bM_n$.
\end{proof}

\section{Invertibles}

\begin{prop}\label{invert}
If $\L \subset \M_n$ is a subspace consisting of singular matrices,
then $\L$ is not transitive.
\end{prop}

\begin{proof} 
Let $k$ be the largest rank of an element of $\L$, 
and fix $A \in \L$ with $\rank A = k$.
Transitivity is unchanged if $\L$ is multiplied on
either side by invertible operators.
So after such a change, we may suppose that
$A=A^2=A^*$ is a projection.

Decompose $\H = A\H \oplus (I-A) \H$.
With this decomposition, each element $L \in \L$ has the form 
$L \simeq \begin{bmatrix} L_0&E\\F&L_1\end{bmatrix}$.
Since $L+tA = \begin{bmatrix} L_0+tI_k &E\\F&L_1\end{bmatrix}$,
the $1,1$ entry is invertible for large $t$, and so it factors as
\[
 \begin{bmatrix} L_0 \!+\! tI_k&0\\F&I_{n-k}\end{bmatrix}
 \begin{bmatrix} I_k &0\\0&L_1 \!-\! F(L_0 \!+\! tI_k)^{-1}E\end{bmatrix}
 \begin{bmatrix} I_k &(L_0 \!+\! tI_k)^{-1}E\\0&I_{n-k}\end{bmatrix} .
\]
From this, it follows that 
\[ k \ge \rank(L+tA) = k + \rank (L_1 - F(L_0+tI_k)^{-1}E).\]
Therefore $L_1 = F(L_0+tI_k)^{-1}E$.
As the right side tends to $0$ as $t\to\infty$,
$L_1=0$ for all $L \in \L$.

This shows that $(I-A) \L (I-A) = 0$ and so $\L$ is not transitive.
\end{proof}

\begin{prop}
If $\L \subset \bM_n$ is transitive, let $r$ be the minimal rank
of non-zero elements of $\L$ and let $s$ be the largest rank of singular elements of $\L$.
Then $r+s \ge n$.
\end{prop}

\begin{proof}
Let $F \in \L$ with $\rank F = r$.
Suppose first that there is an invertible element $A \in \L$ such that
$0 \ne \lambda \in \sigma(A^{-1}F)$.
Then
\[ 0 = \det (\lambda I - A^{-1}F) = \det(A^{-1}) \det (\lambda A - F) .\]
As $ \det(A^{-1}) \ne 0$, $\lambda A - F$ is singular and thus has
rank at most $s$.
But clearly it has rank at least $n-r$.
So $r+s \ge n$.

Otherwise, for every invertible $A$ in $\L$, $A^{-1}F$ is nilpotent.
By our Proposition~\ref{invert}, $\L$ contains invertible elements.
Select an invertible $B \in \L$ so that among the elements
of the form $A^{-1}F$ for invertible $A$ in $\L$, the operator
$F_0 = B^{-1}F$ is nilpotent of the greatest index, say $m+1$.
That is, $F_0^m \ne 0 = F_0^{m+1} = (A^{-1}F)^{m+1}$
for all $A \in \L$ which are invertible.

For any $L\in \L$, and sufficiently small $\mu$, 
the operator $B - \mu L$ is invertible.
So 
\[
 0 = \big((B - \mu L)^{-1}F_0\big)^{m+1} = (I - \mu B^{-1}L)^{-1}F_0 \big((I - \mu B^{-1}L)^{-1}F_0\big)^m
\]
Therefore expanding $(I - \mu B^{-1}L)^{-1} = \sum_{k\ge0} \mu^k (B^{-1}L)^k$,
\[
 0 = F_0 \big((I - \mu B^{-1}L)^{-1}F_0\big)^m = \sum_{k\ge0} \mu^k X_k .
\]
All coefficients of this power series must vanish, and in particular
\[ 0 = X_1 =  F_0B^{-1}LF_0^m + F_0^2B^{-1}LF_0^{m-1} + \dots + F_0^mB^{-1}LF_0 .\]
Multiply on the left by $F_0^{m-1}$ to obtain $F_0^mB^{-1}LF_0^m = 0$
for all $L \in \L$.
This means that $\L \ran F_0^m \subset \ker F_0^m B^{-1}$.  
Both $ \ran F_0^m$ and $\ker F_0^mB^{-1} = B \ker F_0^m$ are proper subspaces,
so this contradicts the transitivity of $\L$.
This contradiction establishes the result.
\end{proof}

\section{Infinite dimensional results}

In this section we examine to what degree the results of the 
previous sections extend to the infinite dimensional setting.   
As we shall see, there are more negative results than positive results.  We begin with an infinite dimensional version of Proposition~\ref{invert}.

Recall that if $T \in \B(\H)$, then an element $\lambda$ of the spectrum 
$\sigma(T)$ of $T$ is called a  \emph{Riesz point} if $\lambda$ is 
an isolated point of $\sigma(T)$, and if 
 $\bigvee_{k\ge1} \ker (\lambda I - T)^k$ is finite dimensional.   
In particular, $\lambda$ is not an element of the essential 
spectrum $\sigma_e(T)$ of $T$.

\begin{thm}\label{invertibles}
Let $\L$ be a subspace of $\B(\H)$ consisting of singular operators.
If $\L$ contains an operator $A$ with $0$ as a Riesz point of its spectrum,
then $\L$ is not topologically transitive.
\end{thm}

\begin{proof}
After a similarity,  $A \simeq \begin{bmatrix} A_0&0\\0&A_1\end{bmatrix}$
where $\H = \M \oplus \M^\perp$ and $\M$ is the spectral subspace
$\ker A^p$ for all sufficiently large $p$.
So $A_1$ is invertible in $\B(\M^\perp)$, and
$A^p \simeq \begin{bmatrix} 0&0\\0&A_1^p\end{bmatrix}$.

Any $L \in \L$ has the form 
$L \simeq \begin{bmatrix} L_0&E\\F&L_1\end{bmatrix}$.
As in the proof of Proposition~\ref{invert}, the $2,2$ entry of $L_t:=L+tA^p$ is invertible
for all $t$ sufficiently large.
Thus $L_t$ is Fredholm of index $0$.
Since $0 \in \sigma(L_t)$, it must have non-trivial kernel.
However, $L_t$  factors as above as
\[
 \begin{bmatrix} I_\M&E(L_1\!+\!tA_1^p)^{-1}\\0&I_{\M^\perp}\end{bmatrix}
 \begin{bmatrix} L_0 \!-\! E(L_1 \!+\! tA_1^p)^{-1}F &0\\0&I_{\M^\perp} \end{bmatrix}
 \begin{bmatrix} I_\M &0\\F&L_1 \!+\! tA_1^p\end{bmatrix} .
\]
The middle factor must have kernel.
Letting $t\to\infty$ shows that $L_0$ is singular.

This shows that $P_\M \L P_\M$ consists of singular matrices.
Hence it is not transitive by Proposition~\ref{invert}.
Since the rank of $P_\M$ is finite,  $\L$ is not topologically transitive, as observed in the last paragraph of the introduction.
\end{proof}

\begin{eg}
The set $\fK$ of compact operators is transitive and singular.

To get even closer to the hypotheses of Theorem~\ref{invertibles},
consider the set $\L = \bC S^* + \fK$, where $S^*$ is the backward shift.
Every element is singular, and there are elements
$A \in \L$ such that $0 \not\in \sigma_e(A)$ and $\ker A \ne 0$.
Nevertheless, this is a transitive space.
\end{eg}


The spaces $\fT$ and $\fH$ of Toeplitz and of Hankel operators 
appear frequently as counterexamples to possible extensions of our finite dimensional results.   
Let $dm$ denote normalized Lebesgue measure on the unit circle 
$\bT = \{ z \in \mathbb{C}: |z| = 1 \}$.   
The set $\{ e_n = e^{i n \theta} \}_{n \in\bZ}$ forms an orthonormal basis 
for the Hilbert space $L^2(\bT) = L^2(\bT, dm)$.   
Let us denote by $\M^\infty (\bT)$ the space 
$\{ M_f: f \in L^\infty(\bT, dm) \}$ of multiplication operators on $L^2(\bT)$.  
It is well-known and  routine to verify that $T = [t_{i j}] \in \M^\infty(\bT)$ 
if and only if $T \in \B(L^2(\bT))$ and $t_{i j} = t_{i + k \ j+k}$ 
for all $i, j, k \in\bZ$.

We denote by $H^2(\bT)$ the Hardy space 
$\overline{\spn} \{ e_n \}_{n=0}^\infty$ of analytic functions in 
$L^2(\bT)$.  
If $P$ denotes the orthogonal projection of $L^2(\bT)$ 
onto $H^2(\bT)$, then the Toeplitz operators are the 
elements of $\fT = P \M^\infty|_{H^2(\bT)}$, 
and with $Q = (I-P)$,  the set of  Hankel operators is 
$\fH = Q \M^\infty|_{H^2(\bT)}$.  


\begin{eg} \label{HankToeptrans}
The spaces $\fH$ and $\fT$ of Hankel and Toeplitz operators  are topologically transitive.

\bigskip

We show that the space $\M^\infty P$ is topologically transitive.   Since $\fT$ and $\fH$ are compressions of $\M^\infty P$, it immediately follows that they too are topologically transitive.

Consider $A \in (\M^\infty P)_\perp$ with $\rank A = 1$.   Then $A = f g^*$ for some $f \in H^2(\bT)$ and $g \in L^2(\bT)$.   Since $A \ne 0$, neither $f$ nor $g$ is zero.   The condition $A \in (\M^\infty P)_\perp$ implies that $\mathrm{tr} (M_h P A) = \mathrm{tr} (M_h f g^*) = \langle h f, g \rangle  = 0$ for all $h \in L^\infty(\bT)$.   That is, $\int_\bT h f \overline{g} = 0$ for all $h \in L^\infty(\bT)$.   Since $f \overline{g} \in L^1(\bT) = (L^\infty (\bT))_\perp$, it follows that $f \overline{g} = 0$ a.e..   The classical F.~and~M.~Riesz Theorem (see \cite[Theorem~6.13]{Douglas}) asserts that as $0 \ne f \in H^2(\bT)$, the set $\{ z \in \bT: f(z) = 0 \}$ has measure $0$.   From this it follows that $g = 0$ a.e., a contradiciton.   Thus $\M^\infty P$ is topologically transitive.

We observe that an analogous argument may be used to establish the fact that $Q \M^\infty$ is also topologically transitive.

Note also that $\fT e_0 = \{ T_h  e_0 = P h: h \in L^\infty(\bT) \}$ is dense in $H^2(\bT)$, but is not everything.   Thus $\fT$ is an example of a topologically transitive space which is not transitive.
\end{eg}


\begin{defn} \label{totallysep}
A subspace $\L \subseteq \B(\H_1, \H_2)$ is said to be  \emph{totally separating} 
if $\L$ is $k$-separating for all $k \geq 2$.
\end{defn}

There is no point in defining \textit{totally topologically transitive} in the analogous way,
because this would just say that $\L$ is \sot-dense in $\B(\H)$.
A natural modification of Example~\ref{nsepnottrans} shows that for each $k \geq 1$ there are
totally separating spaces which are topologically $k$-transitive but not topologically $(k+1)$-transitive.


\begin{prop} \label{prod_trans2sep}
Let $\L_1$ and $\L_2$ be subspaces of $\B(\H)$, and suppose that $\L_1$ is topologically transitive.
\begin{enumerate}
	\item[(a)]
	If $\L_2$ is $k$-separating for some $k \ge 1$, then $\spn \L_1 \L_2$ is topologically $k$-transitive.
	\item[(b)]
	If $\L_2$ is totally separating, then $\spn \L_1 \L_2$ is dense in the strong operator topology on $\B(\H)$.
\end{enumerate}	
\end{prop}

\begin{proof}
(a) Choose $x_1, x_2, ..., x_k \in \H$ linearly independent.   
Let $\ep > 0$, and choose $y_1, y_2, ..., y_k \in \H$ arbitrary.
Since $\L_2$ is $k$-separating, we can find operators $L_1, L_2, ..., L_k \in \L_2$ 
so that $L_i x_i \not = 0$, but $L_i x_j = 0$ for all $1 \leq i \not = j \leq k$.  
Since $\L_1$ is topologically transitive, we can find $K_1, K_2, ..., K_k \in \L_1$ with 
$\| K_j (L_j x_j) - y_j \| < \ep/k$, $1 \leq j \leq k$.

Let $T = \sum_{j=1}^k K_j L_j \in \L_1 \L_2$.   Then 
\[
\| T x_j - y_j \| = \| ( \sum_{i=1}^n  K_i L_i ) x_j - y_j \| = \| K_j L_j x_j - y_j \| < \ep/k, \]
for each $1 \leq j \leq k$.  Since $\{ x_j \}_{j=1}^k$ linearly independent and 
$\{y_j\}_{j=1}^k$ are arbitrary, $\spn \L_1 \L_2$ is topologically $k$-transitive.

(b)
By part (a), $\spn \L_1 \L_2$ is topologically $k$-transitive for all $k \ge 1$.  
By the comments preceding the proposition, this says that
$\spn \L_1 \L_2$ is dense in the strong operator topology in $\B(\H)$.
\end{proof}




\begin{eg}  \label{Hanknot2trans}
The Hankel operators $\fH$ and the Toeplitz operators $\fT$ are
products of two topologically $1$-transitive spaces,
but they are not topologically $2$-transitive nor even 2-separating.

Indeed, $\fH = (Q \M^\infty) (\M^\infty P)$  and $\fT = (P \M^\infty) (\M^\infty P)$.   
In finite dimensions, we have seen that the product of two transitive spaces is $2$-transitive.  
We have also seen that in the finite dimensional setting, $1$-transitive spaces 
are automatically $2$-separating.   

A typical operator in $\fH$ admits an infinite matrix representation of the form
\[
H = \begin{bmatrix} a_1 & a_2 & a_3 & ... \\ a_2 & a_3 & ... & \\
 a_3 & ... & ... & ... \\ \vdots & & & \end{bmatrix} 
\]
relative to the bases $\{ e_n \}_{n=0}^\infty$ for $H^2(\bT)$ and $\{ e_n \}_{n = -1}^\infty$ for $(H^2(\bT))^\perp$.
Therefore $H e_0 = 0$ implies $H = 0$, 
and hence $\fH $ is neither $2$-separating nor topologically $2$-transitive.

Similarly, if $T\in\fT$ and $Te_n=0$, then $Te_i=0$ for $0 \le i < n$.
So $\fT$ is also neither 2-separating nor topologically $2$-transitive.
This can be contrasted with Example~\ref{Toeplitz2sepnot3} where 
it is shown that $\fT_n$ is 2-separating but not 3-separating.
\end{eg}



The following technical result will be used in the proof of Proposition~\ref{prodtransHank}.

\smallskip

\begin{lem} \label{manyranks}
Suppose $\L_2 \subseteq \B(\H)$ contains a sequence $(F_m)_{m=1}^\infty$ of operators with 
\begin{enumerate}
	\item[1)]
	$\mathrm{rank\,} F_m = m$, $m \geq 1$;
	\item[2)]
	$\ker F_{m+1} \subseteq \ker F_m$ for all $m \geq 1$; and 
	\item[3)]
	$\bigcap_{m \geq 1} \mathrm{ker\,} F_m = \{ 0 \}$.
\end{enumerate}
If $\L_1$ is topologically transitive, then $\spn \L_1 \L_2$ is dense in the weak operator topology on $\B(\H)$.
\end{lem}
	
\begin{proof}
Choose a sequence $(F_m)_{m=1}^\infty \subseteq \L_2$ as in the statement of the Lemma.
For each $m \geq 1$, let $\H_m = (\ker F_m)^\perp$, so that $\H_m \subseteq \H_{m+1}$, 
and $\mathrm{dim\,} \H_m = m$ for all $m$.  
Fix $e_1 \in \H_1$ with $\| e_1 \| = 1$, and for $m \geq 2$, 
choose $e_m \in \H_m \ominus H_{m-1}$ with $\| e_{m} \| = 1$.   
The third hypothesis above guarantees that $\{ e_m \}_{m=1}^\infty$ 
spans $\H$, and thus forms an orthonormal basis for $\H$.  
Our goal is to show that if $P_m$ is the orthogonal projection of $\H$ onto $\H_m$,
then $\overline{\spn \L_1 \L_2}^\wot$ contains $\B(\H) P_m$ for all $m \geq 1$.   
Since the latter set is clearly dense in the weak operator topology, so is the former.

Let $T \in \B(\H)$ be arbitrary, and set $\ep > 0$.   
Now $e_m \in \H_m$, so $z_m := F_m e_m \not = 0$.   
Fix $R_m \in \L_1$ so that $\| R_m F_m e_m - T e_m \| < \ep/m$.  
Next, $e_{m-1} \in \H_{m-1}$, so $F_{m-1} e_{m-1} \not = 0$.  
Fix $R_{m-1} \in \L$ so that 
\[
 \| R_{m-1} z_{m-1} - (T e_{m-1} - R_m F_m e_{m-1}) \| < \ep/m. 
\]
Observe that $R_{m-1} F_{m-1} e_m = 0$, so that 
\[ \| ( R_{m-1} F_{m-1} + R_m F_m) e_m - T e_m \| < \ep/m .\] 
More generally, having chosen $R_m, R_{m-1}, ..., R_{m-k}$, 
we can choose\linebreak $R_{m-(k+1)}$ in $\L_1$ so that 
\[
 \big\| R_{m-(k+1)} F_{m-(k+1)} e_{m-(k+1)}
 -\big( T - \sum_{j=m-k}^m R_j F_j \big) e_{m-(k+1)} \big\| < \ep/m. 
\]
 	
It follows  that $Q_m = \sum_{j=1}^m R_j F_j$ satisfies
$\| Q_m e_r - T e_r \| < \ep/m$ for $1 \leq r \leq m.$
Since $Q_m = Q_m P_m$,  $\| Q_m - T P_m \| <  \ep$.
Finally, since $Q_m$ belongs to $\spn \L_1\L_2$ \vspace{.3ex} and $\ep>0$ is arbitrary,
$\B(\H) P_m \subset \ol{\spn \L_1\L_2}$.
Therefore $\overline{\spn \L_1 \L_2}^\wot = \B(\H)$.
\end{proof}
	

\begin{prop} \label{prodtransHank}
Suppose that $\L \subseteq \B(\H)$ is topologically transitive,  
and let $\fH \subseteq \B(\H)$ denote the space of Hankel operators.
Then 
\[ 
\overline{ \spn \L \fH}^\wot = \overline{ \spn  \fH \L }^\wot = \B(\H) .
\]  
In particular, $\spn \fH^2$ is weak operator dense in $\B(\H)$.
\end{prop}

\begin{proof}
For each $n \geq 1$, the rank $n$ operators 
\[
F_n = \begin{bmatrix} 1 & 1 & 1 & ... & 1 & 0 & 0 & ... \\
				   1 &1 & ... & 1 & 0 & 0 & ... & ... \\
				   \vdots & & & & & & & ... \\	
				   1  & 1 & ... & & & & & \\
				   1 & 0 & & & & & & \\
				   0 & ... & & & & & & \\
	\end{bmatrix}
	\]
lie in $\fH$ and satisfy the conditions of Lemma~\ref{manyranks}.   
Thus we may conclude that $\overline{ \spn \L \fH}^\wot = \B(\H)$.   

Now $\L$ topologically transitive implies that $\L^\mathrm{t}$ is topologically transitive.   
Since $ \fH = \fH^\mathrm{t}$, it follows that 
$(\fH \L)^\mathrm{t} = \L^\mathrm{t} \fH$ 
is weak operator dense in $\B(\H)$, whence $\overline{ \spn  \fH \L }^\wot = \B(\H)$. 
\end{proof}			   

\bigskip

Let us next consider $\spn \fT^2$.   Let $E_{ij} = e_i e_j^*, \ i, j \ge 0$ denote the matrix units of $\B(H^2(\bT))$.   Letting $S = P M_z|_{H^2(\bT)} \in \fT$, $S$ is unitarily equivalent to the  unilateral forward shift, and a routine calculation reveals that for $i , j \ge 0$, 
\begin{align*}
 E_{i, j} = S^i (I - S S^*) (S^*)^j = S^i (S^*)^j - S^{i+1} (S^*)^{j+1}. 
\end{align*}

Since $S^{k}, (S^*)^{l} \in \fT$ for all $k, l \ge 0$, $E_{i j} \in \spn \fT^2$ for all $i, j \ge 0$.  Thus the norm closure of $\spn \fT^2$ contains all compact operators, and is therefore transitive.

So far we have not been able to determine whether $\L \subseteq \B(H^2(\bT))$ topologically transitive implies that $\overline{ \spn  \L \fT }^\wot = \B(\H)$.

\bigskip

As we have seen in Section~\ref{section4}, if $\L \subseteq \bM_n$ is transitive, then $\spn \L^r = \bM_n$ for some $1 \le r \le n$.  It would be interesting to find estimates for $\kappa_n := \min \{ 1 \le r \le n:   \spn \L^r = \bM_n \mbox{ for all } \L \subseteq \bM_n \mbox{ transitive}\}$.  In particular, is $( \kappa_n)_{n=1}^\infty$ bounded?    In the infinite dimensional setting, does there always exist some $r \ge 1$ so that $\L \subseteq \B(\H)$ topologically transitive implies that $\overline{ \spn \L^r }^\wot = \B(\H)$?   More generally, does there exist $s \ge 1$ so that if $\L_1, \L_2, ..., \L_s$ are topologically transitive subspaces of $\B(\H)$, then $\overline{ \spn  \L_1 \L_2 \cdots \L_s }^\wot = \B(\H)$?    By Example~\ref{Hanknot2trans}, if such an $s$ exists, then $s \ge 3$.

\begin{eg} \label{n-1-transnotnsep} A subspace of $\B(\H)$ can be topologically $(n-1)$-transitive, 
but not $n$-sepa\-rating.

Let $\{e_n = z^n: n \in \bZ\}$ be the standard orthonormal basis for $L^2(\bT)$,  and set $\W_n = \spn \{e_1, e_2, ..., e_n \}$.   
Consider the space $\L_n  \subseteq \B(\W_n, L^2(\bT))$ of operators 
of the form  $A = [ a_{ij}]$ where $\sum_{k=1}^n a_{i+k,k} = 0$, $i \in\bZ$.
To see that $\L_n$ is not $n$-separating, it suffices to observe that if 
the first $(n-1)$ columns of $\L_n$ are zero, then the last column of 
$\L_n$ is necessarily zero as well.

The proof that $\L_n$ is topologically $(n-1)$-transitive relies upon the structure of $(\L_n)_\perp$.  
We may, in a manner analogous to that used in our analysis of the finite dimensional setting, 
identify $(\L_n)_\perp$ with the set of trace class 
operators $B = [b_{i j}] \subseteq \B(\W_n, L^2(\bT))$ which 
satisfy $\sum_{i=1}^n \sum_{j \in\bZ} a_{i j} b_{ji } = 0$ for all $A = [a_{i j}] \in \L_n$.
A routine calculation then shows that $B \in (\L_n)_\perp$ implies that 
$b_{k, i + k} = b_{1, i+ 1}$ for all $2 \leq k \leq n$.  

If we think of the rows of $B$  as vectors in $L^2(\bT)$, 
then this says that $B \in (\L_n)_\perp$ if and only if 
$B^t = \begin{bmatrix} f & z f & z^2 f & ... & z^{n-1} f \end{bmatrix}$ 
for some $f \in L^2(\bT)$.   
If $\L_n$ were not topologically $(n-1)$-transitive, then we could find such a 
$B\ne0$ with $\rank B \leq n-1$.  
In particular, $\rank B^t = \rank B \leq n-1$ and so $\ker  B^t \ne 0$.   
Choose a vector $0 \ne \sum_{i=1}^n \alpha_i e_i \in \ker B^t$.
Then $(\sum_{i=1}^n \alpha_i z^{i-1})  f = 0$ a.e.
Since $\sum_{i=1}^n \alpha_i z^{i-1}$ is a non-trivial polynomial, 
it has at most finitely many zeroes.
Therefore $f = 0$ a.e., contradicting $B \ne 0$.   
Hence $\L_n$ is topologically $(n-1)$-transitive.
\end{eg}


\begin{eg} \label{intersecttrans}
The intersection of a descending sequence of \wot-closed transitive spaces 
need not be topologically transitive.

As before, we let $\{e_n : n \in \bZ\}$ be the standard basis for $L^2(\bT)$.
Let $R_n$ denote the orthogonal projection of $L^2(\bT)$ onto $\spn \{ e_k : |k|\le n \}$.  
Then $R_n \M^\infty|_{R_n L^2(\bT)} \subseteq \B(R_n L^2(\bT))$ is clearly 
unitarily equivalent to the Toeplitz matrices on $\H_{2n+1}$, 
and so it is transitive (see Example~\ref{toeplitz_are_trans}).

Let $\R_n = \{ X \in \B(L^2(\bT)): R_n X R_n \in R_n \M^\infty R_n \}$.   
Then $\R_n$ is transitive and \wot-closed. 
Indeed, with respect to the decomposition 
$L^2(\bT) = R_n L^2(\bT) \oplus (R_n L^2(\bT))^\perp$, 
an element of $\R_n$ looks like $\begin{bmatrix} X_1 & X_2 \\ X_3 & X_4 \end{bmatrix}$
where $X_1 \in R_n \M^\infty|_{R_n L^2(\bT)}$ and the other entries are arbitrary.
As the matrix entries are independent and each corner is transitive,
it follows easily that $\R_n$ is transitive.

Observe, however, that $\bigcap_{n=1}^\infty \R_n = \M^\infty(\bT)$.   
Since $\M^\infty(\bT)$ has many proper closed invariant subspaces, it is not topologically transitive.  So 
the intersection of a descending sequence of transitive spaces need not be topologically transitive.

Note that there are limits to the decreasing intersection of 
transitive spaces $\L = \bigcap_{n\ge1}\L_n$.
For if $P$ and $Q$ are rank $n$ projections, then $P\L_n Q$ is a transitive
subspace of $\B(Q\H,P\H)$ for all $n$, and so has dimension at least $2n-1$.
Thus the same is true for the intersection.
Moreover since a decreasing sequence of subspaces of a finite dimensional
space is eventually constant, we see that $P\L Q$ is transitive whenever
$P$ and $Q$ are finite rank.
Our example shows that this estimate is sharp because the compression
using $P\H=Q\H = \spn\{e_i:0\le i < n\}$ yields $\fT_n$ as the intersection;
and it has dimension exactly $2n-1$.
\end{eg}

\begin{thm}\label{T:infinite_product}
Suppose that $\L, \M \subset\B(\H)$ are topologically transitive.  
If $\L$ is contained in the \wot-closed span of its rank one elements, 
then the norm closure of $\spn \L \M$ is transitive.
\end{thm}

\begin{proof}
For each rank one element $uv^* \in \L$, $\L\M$ contains $uv^*M = u (M^*v)^*$
for all $M \in \M$.  
By the topological transitivity of $\M$, the norm closure of $\L \M$ contains $ue_j^*$ where $\{e_j\}$ is an orthonormal basis for $\H$.
As $\L$ is topologically transitive and \wot-spanned by rank ones, the collection of such vectors $u$ densely spans $\H$, from which it follows that the norm closure of  $\spn \L\M$ contains the compact operators.
\end{proof}


\begin{thm}\label{T:infinite_tensor}
Suppose that $\L \subset \B(\H_1,\H_2)$ and $\M \subset\B(\K_1,\K_2)$
are topologically transitive.  If $\L$ is contained in the \wot-closed span of
its rank one elements, then the spatial tensor product $\L\otimes\M$ is topologically transitive.
\end{thm}

\begin{proof}
Fix $0 \ne x \in \H_1\otimes \K_1$ and $y \in \H_2\otimes \K_2$.
If $L=uv^*$ is a rank one element of $\L$, then $(L \otimes I_{\K_1}) x = u\otimes z_L$ for some $z_L \in \K_1$.
We claim that 
\[ \spn \{u \in \H_2 : uv^* \in \L \AND (uv^*\otimes I_{\K_1}) x \ne 0 \} = \H_2 .\]
Indeed, if $0 \ne w\in\H_2$ is a vector orthogonal to this span, 
then for all $z \in \K_1$ and all rank one $L \in \L$,
\[ 0 = \ip{(L \otimes I_{\K_1}) x, w \otimes z} = \ip{x, (L^* \otimes I_{\K_1})(w \otimes z)} 
     = \ip{x, L^* w \otimes z} .\]
Since the span of the  rank one elements is \wot-dense in $\L$, we deduce
that $\spn \{L^*w : L \in \L,\ \rank L = 1 \} = \ol{\L^* w} = \H_1$.
Since $z$ is arbitrary in $\K_1$, it follows that the vectors of the form
$L^* w \otimes z$ span $\H_1\otimes\K_1$, implying that $x=0$.

Approximate $y$ by a finite sum $\sum_{i=1}^N u_i \otimes y_i$
where $L_i = u_i v_i^* \in \L$ and $(u_i v_i^* \otimes I_{\K_1}) x = u_i \otimes z_i \ne 0$.
Choose $M_i \in \M$ so that $M_i z_i \approx y_i$.
Then $\L \otimes\M$ contains $A = \sum_{i=1}^N L_i \otimes M_i$ and
\[
 A x = \sum_{i=1}^N (I_{\H_2} \otimes M_i)(L_i\otimes I_{\K_1}) x
 = \sum_{i=1}^N (I_{\H_2} \otimes M_i) u_i \otimes z_i \approx \sum_{i=1}^N u_i \otimes y_i .
\]
Take appropriate limits to complete the proof.
\end{proof}


If the subspace has additional structure, such as being a module over a masa, 
then stronger results may hold.  
For example, a non-trivial result of Arveson~\cite{Arv} 
(see also \cite[Theorem~15.9]{Dnest}) shows:

\begin{thm}\label{bimodule}
Let $\D_i$ be masas in $\B(\H_i)$.
A topologically transitive subspace $\L$ of $\B(\H_1,\H_2)$ which is a
$\D_2$--$\D_1$ bimodule is  \wot-dense in $\B(\H_1,\H_2)$.
\end{thm} 

We provide a new, more elementary proof.
Actually Arveson's proof works for the weak-$*$ topology,
whereas this proof is only valid for the \wot-topology.

\smallskip

\begin{lem}~\label{rectangle}
Let $(X_i,\mu_i)$ be regular Borel measures.
Let $k(x,y) = \sum_{i=1}^m \alpha_i(x) \beta_i(y) \in L^2(\mu_1 \times \mu_2)$
where $\alpha_i \in L^2(\mu_1)$ and $\beta_i \in L^2(\mu_2)$.
For any $a$ in the essential range of $k$ and any $\ep>0$,
there is a measurable rectangle $A_1 \times A_2$ with
$0 < \mu_i(A_i) < \infty$ such that $|k(x,y)-a|<\ep$ for all
$(x,y) \in A_1 \times A_2$.
\end{lem}

\begin{proof}
Choose a measurable rectangle $Y_1 \times Y_2$ of finite positive measure
on which $k$ is bounded and still has $a$ in its essential range.
It is a standard argument to approximate each $\alpha_i$ and $\beta_i$
uniformly (and in $L^2$) by simple functions on $Y_1$ and $Y_2$ respectively.
Combining these simple functions allows us to approximate $k\upchi_{Y_1\times Y_2}$ uniformly by 
a finite linear combination of characteristic functions of measurable rectangles.
We may then pick a rectangle on which $k$ takes values close to $a$.
\end{proof}

\begin{rem}
The lemma fails for arbitrary functions in $L^2(\mu_1 \times \mu_2)$.
For example. take $\mu_1=\mu_2$ to be Lebesgue measure on $[0,1]$.
Let $A$ be a compact nowhere dense subset of $[0,1]$ with positive measure.
Then $k(x,y) = \upchi_A(x-y)$ has $1$ in its essential range.
However if $k = 1$ on a measurable rectangle $A_1 \times A_2$,
then $A_1 - A_2 \subset A$ is nowhere dense.
It is a well known fact that the difference of two measurable sets of
positive measure has interior.  So $A_1 \times A_2$ has measure $0$.
\end{rem}

\begin{proof}[\em\textbf{Proof of Theorem~\ref{bimodule}.}]
By the spectral theorem for masas, we may suppose that there are regular Borel
spaces $(X_i,\mu_i)$ so that $\D_i$ are unitarily equivalent to $L^\infty(\mu_i)$
acting by multiplication on $\H_i = L^2(\mu_i)$.
If $\L$ is not \wot-dense, then there is a finite rank operator $F \in \L_\perp$.
Moreover, it is evident that $\L_\perp$ is a $\D_1$--$\D_2$ bimodule.
Our goal is to show that using $F$ and the bimodule property,
we may find a rank one element of $\L_\perp$.
This will contradict topological transitivity.

Observe that $F$ may be written as an integral operator with kernel 
$k(x,y) = \sum_{i=1}^m \alpha_i(x) \beta_i(y)$
where $\alpha_i \in L^2(\mu_1)$ and $\beta_i \in L^2(\mu_2)$.
Since $F \ne 0$, there is a non-zero value $a$ in the essential range of $k$.
By Lemma~\ref{rectangle}, there is a measurable rectangle $A_1 \times A_2$ 
of finite non-zero measure so that $|k(x,y)-a|<|a|/2$ for all
$(x,y) \in A_1 \times A_2$.

Let $h(x,y) = \upchi_{A_1 \times A_2} k(x,y)^{-1}$.
Then $h \in L^\infty(\mu_1 \times \mu_2)$.
Hence $h$ is a limit in $L^\infty(\mu_1 \times \mu_2)$ of a sequence
of simple functions of the form $h_k = \sum_{j=1}^{m_k} f_{kj}(x) g_{kj}(y)$.
It follows by routine calculations that
\[
 \sum_{j=1}^{m_k} M_{f_{kj}\chi_{A_1}} F M_{g_{kj}\chi_{A_2}}
 \quad\text{has kernel}\quad
 \sum_{j=1}^{m_k} f_{kj}(x) k \upchi_{A_1\times A_2} g_{kj}(y) .
\]
This converges in $L^2(\mu_1\times\mu_2)$ to $\upchi_{A_1\times A_2}$,
and thus the corresponding operators converge in norm to the rank one
integral operator with kernel $\upchi_{A_1\times A_2}$.
This produces a rank one element of $\L_\perp$.
\end{proof}

The following result is very easy.
Recall that a masa is atomic if it consist of all diagonal operators with
respect to some orthonormal basis.

\begin{prop} \label{atomicmasamodule}
Suppose that a topologically transitive subspace $\L \subseteq \B(\H)$ is a 
left or right module over an atomic masa $\fD$. Then $\overline{\L}^\wot = \B(\H)$.
\end{prop}

\begin{proof}
First we suppose that $\L$ is a right $\fD$-module.
Let $\fD$ be diagonal with respect to $\{e_n : n\ge1\}$;
and let $P_n = e_ne_n^*$.
Since $\L$ is topologically transitive, $\ol{\L P_n} = \B(\H)P_n$.
Summing over $n$ yields a \wot-dense subspace.

By considering $\L^t$ which is also topologically transitive, 
we obtain the left $\fD$-module case.
\end{proof}


\begin{eg}\label{cntsmasamodule}
A right or left module over a non-atomic masa need not be \wot-dense.

Example~\ref{HankToeptrans} showed that $\M^\infty P \subset\B(H^2,L^2(\bT))$
is topologically transitive.  
It is evidently a left $\M^\infty$ module, and a proper \wot-closed subspace.
The adjoint $P\M^\infty \subset \B(L^2(\bT), H^2)$ is similarly a right
 $\M^\infty$ module, and a proper \wot-closed subspace.
\end{eg}


\begin{prop}  \label{prodmasamodules}
Suppose that $\L$ is a left module \vspace{.3ex} over a masa $\D$,  
and that $\L$ is topologically transitive.  
Then $\overline{\spn \L \fD}^\wot = \ol{ \spn \L^2}^\wot = \B(\H)$.
\end{prop}

\begin{proof}
Evidently $\spn \L \fD$ is a $\fD$-bimodule.  So it is \wot-dense by Arveson's Theorem.
Therefore 
\begin{align*}
 \ol{ \spn \L^2}^\wot &= \ol{\spn \L\fD \L}^\wot = \ol{\spn \B(\H) \L}^\wot = \B(\H)
 .\qquad\quad \qedhere
\end{align*}
\end{proof}


\end{document}